\documentclass{amsart}
\usepackage{times,amssymb,array}

\def\Q{\mathbb Q}

\def\Z{\mathbb Z}
\def\SL{\mathrm{SL}}
\def\tm{{\tilde m}}
\def\tn{{\tilde n}}
\def\gen#1{\langle#1\rangle}
\def\Mat#1#2#3#4{\begin{pmatrix}#1&#2\\#3&#4\end{pmatrix}}
\def\SMat#1#2#3#4{\left(\begin{smallmatrix}#1&#2\\#3&#4\end{smallmatrix}
\right)}
\def\Re{\mathrm{Re}}

\newtheorem{Lemma}{Lemma}[section]

\newtheorem{Theorem}{Theorem}
\theoremstyle{definition}
\newtheorem{Definition}[Lemma]{Definition}
\newtheorem{Notation}[Lemma]{Notation}
\newtheorem{Remark}[Lemma]{Remark}

\begin{document}

\title{Twisted Hecke $L$-values and period polynomials}

\author{Shinji Fukuhara}
\address{Department of Mathematics, Tsuda College, Tsuda-machi 2-1-1,
  Kodaira-shi, Tokyo 187-8577, Japan}
\email{fukuhara@tsuda.ac.jp}
\author{Yifan Yang}
\address{Department of Applied Mathematics, National Chiao Tung University,
  1001 Ta Hsueh Road, Hsinchu, Taiwan 300}
\email{yfyang@math.nctu.edu.tw}
\date{\today}

\thanks{}

\begin{abstract} Let $f_1,\ldots,f_d$ be an orthogonal basis for the
  space of cusp forms of even weight $2k$ on $\Gamma_0(N)$. Let
  $L(f_i,s)$ and $L(f_i,\chi,s)$ denote the $L$-function of $f_i$ and
  its twist by a Dirichlet character $\chi$, respectively. In this
  note, we obtain a ``trace formula'' for the values
  $L(f_i,\chi,m)\overline{L(f_i,n)}$ at integers $m$ and $n$ with $0<m,n<2k$ and
  proper parity. In the case $N=1$ or $N=2$, the formula gives us a
  convenient way to evaluate precisly the value of the ratio
  $L(f,\chi,m)/L(f,n)$ for a Hecke eigenform $f$.
 \end{abstract}

 \maketitle

 \begin{section}{Introduction and statements of results} Let
   $f(z)=\sum_{n=1}^\infty a_f(n)e^{2\pi inz}$ be a Hecke eigenform of
   even weight $2k$ on $\Gamma_0(N)$ and $f_\chi(z)=\sum_{n=1}^\infty
   a_f(n)\chi(n)e^{2\pi inz}$ be its twist by a Dirichlet character
   $\chi$. The $L$-function, defined by $L(f,s)=\sum_{n=1}^\infty
   a_f(n)n^{-s}$ and extended analytically to the whole complex plane,
   and its twist $L(f,\chi,s)=\sum_{n=1}^\infty a_f(n)\chi(n)n^{-s}$
   are very important number theoretical objects. For instance, when
   $f(z)$ is the weight $2$ newform associated to a rational elliptic
   curve $E$, the Birch and Swinnerton-Dyer conjecture asserts that the
   rank of the group of rational points $E(\Q)$ on $E$ is equal to the
   order of $L(f,s)$ at $s=1$.

   In this article, we are concerned with the values of $L(f,s)$ and
   $L(f,\chi,s)$ at integers inside the critical strip $0<\Re\, s<2k$.
   In \cite{Manin2}, Manin showed that for a normalized Hecke eigenform
   $f(z)$ of weight $2k$ on $\SL(2,\Z)$, there are two real numbers
   $\omega_f^+$ and $\omega_f^-$, depending only on $f$, such that
   $$
     \pi^{-2n}L(f,2n)/\omega_f^+, \qquad
     \pi^{-(2n-1)}L(f,2n-1)/\omega_f^-
   $$
   are contained in the (totally real) field $\Q(a_f(2),a_f(3),\ldots)$
   for all integers $n$ with $0<2n,2n-1<2k$. This result was later
   generalized to newforms on $\Gamma_0(N)$ by Razar
   \cite[Theorem 1]{Razar}. More generally, it is known that the
   twisted Hecke $L$-value $L(f,\chi,n)$ is equal to an algebraic
   number times either $\pi^n\omega_f^+$ or $\pi^n\omega_f^-$,
   depending on the parities of $n$ and $\chi$.

   The values of $L(f,s)$ and $L(f,\chi,s)$ at the center point $s=k$
   are particularly interesting. Assume that $g(z)=\sum_{n=1}^\infty
   b_g(n)e^{2\pi inz}$ is the modular form of weight $k+1/2$ lying in
   Kohnen's plus-space corresponds to $f(z)$ in the sense of Shimura.
   In \cite{Waldspurger}, Waldspurger proved that $b_g(n)^2$ is
   essentially proportional to the value of $L(f,\chi_{(-1)^kn},s)$ at
   $s=k$, where $\chi_D=\left(\frac D\cdot\right)$. Later on, Kohnen
   and Zagier \cite{KZ1} made this result more explicitly by proving
   $$
     \frac{b_g(n)^2}{\gen{g,g}}
    =\frac{n^{k-1/2}\Gamma(k)}{\pi^k}\frac{L(f,\chi_{(-1)^kn},k)}{\gen{f,f}}
   $$
   for a normalized Hecke eigenform $f(z)$ on the full modular group
   $\SL(2,\Z)$ and positive integers $n$ such that $(-1)^kn$ is a
   fundamental discriminant, where $\gen{f,f}$ and $\gen{g,g}$ denote
   the Petersson norms of $f(z)$ and $g(z)$, respectively. This result
   was generalized by \cite{Baruch-Mao} and \cite{Mao} to Hecke
   eigenforms on $\Gamma_0(N)$.

   In this article, we will derive a ``trace formula''
   $$
     \sum_{i=1}^s\frac1{\gen{f_i,f_i}}L(f_i,\chi,m)\overline{L(f_i,n)}
   $$
   for a Dirichlet character $\chi$ and integers $m$ and $n$ with
   proper parity, where $\{f_1,\ldots,f_s\}$ is any orthogonal basis for
   $S_{2k}(\Gamma_0(N))$. In some cases,
   such as $N=1$ and $N=2$, this formula enables us to compute the
   exact value of the ratio $L(f,\chi,m)/L(f,m)$ for a Hecke eigenform
   $f$.

   To achieve our goal, we first express the values of a Hecke
   $L$-function $L(f,s)$ as \emph{periods}
   $$
     r_n(f):=\int_0^{i\infty} f(z)z^n\,dz=\frac{n!}{(-2\pi i)^{n+1}}
     L(f,n+1)
   $$
   of a cusp form $f$. The periods are studied extensively in
   \cite{E1,Fukuhara-AA1997,Fukuhara-MAnn1998,FY,Haberland,KZ2,Manin1,Manin2,SH1,Z1,Z2}. In particular, for the case
   of $\SL(2,\Z)$, Kohnen and Zagier \cite{KZ2} showed that there is
   a rational structure associated to the periods that is different
   from the usual rational structure coming from the Fourier
   coefficients of cusp forms. The idea is to consider the cusp form
   characterized by the property
   $$
     r_n(f)=\gen{f,R_n}
   $$
   for all $f\in S_{2k}(\Gamma_0(N))$, where $\gen{\cdot,\cdot}$ denote
   the Petersson inner product. Then in \cite{FY,KZ2} it is shown that
   the values of $r_m(R_n)$ can be expressed in terms of the Bernoulli
   numbers for integers $m$ and $n$ with opposite parity satisfying
   $0\le m\le 2k-2$ and $1\le n\le 2k-3$. In \cite{F5}, by considering
   the natural correspondence of $S_{2k}(\SL(2,\Z))$, its dual, and the
   space of Dedekind symbols, the first author of the present article
   found bases for $S_{2k}(\SL(2,\Z))$ in terms of $R_n$, which in turn
   give explicit expression for Hecke operators in terms of Bernoulli
   numbers and sum-of-divisor functions. For the case $\Gamma_0(2)$,
   this is done in \cite{FY} with a different approach.

   Now if $\{f_1,\ldots,f_s\}$ is an orthogonal basis for
   $S_{2k}(\Gamma_0(N))$, then we have
   $$
     R_n=\sum_{i=1}^s\frac{\gen{R_n,f_i}}{\gen{f_i,f_i}}f_i,
   $$
   from which we deduce that
   $$
     r_m(R_n)=C_{m,n}\sum_{i=1}^s\frac1{\gen{f_i,f_i}}
     L(f_i,m+1)\overline{L(f_i,n+1)}
   $$
   for some complex number $C_{m,n}$ depending only on $m$ and $n$. In
   other words, the ``trace'' of
   $L(f_i,m+1)\overline{L(f_i,n+1)}/\gen{f_i,f_i}$ is essentially
   $r_m(R_n)$. More generally, if we define the \emph{twisted period}
   of a cusp form $f$ by
   $$
     r_{m,\chi}(f):=\int_0^{i\infty}f_\chi(z)z^m\,dz,
   $$
   where $f_\chi$ denotes the twist of $f$ by $\chi$, then the trace of
   $L(f_i,\chi,m+1)\overline{L(f_i,n+1)}$ is essentially
   $r_{m,\chi}(R_n)$.

   It turns out that formulas for $r_{m,\chi}(R_n)$ can be more
   elegantly stated if we write them collectively as \emph{twisted
   period polynomials}
   $$
     r_{\chi}(f)(X):=\int_0^{i\infty}f_\chi(z)(X-z)^{2k-2}\,dz
   $$
   of a cusp form $f$. Before we state our formula for
   $r_\chi(f)(X)$, let us first fix some notations.

 \begin{Notation} Throughout the notes, the letter $N$ will always
   denote the level of the congruence subgroup $\Gamma_0(N)$, and
   $\chi$ will represent a primitive Dirichlet character modulo $D$
   with $D>1$.

   For convenience, we shall write the weight of the space of cusp
   forms $S_{2k}(\Gamma_0(N))$ under consideration as $2k=w+2$. For
   integers $m$ and $n$, we set
   $$
     \tm=w-m, \quad \tn=w-n.
   $$
 \end{Notation}

   We now recall some definitions related to a Dirichlet character
   $\chi$.

 \begin{Definition} For a non-negative integer $k$, the $k$th Bernoulli
   polynomial $B_k(x)$ is defined by the power series expansion
 $$
   \frac{te^{xt}}{e^t-1}=\sum_{k=0}^\infty\frac{B_k(x)}{k!}t^k.
 $$
 Naturally, $B_k(x)$ is the zero polynomial if $k$ is a negative
 integer. For a Dirichlet character $\chi$ modulo $D$, we also define
 \emph{generalized Bernoulli polynomial} $B_{k,\chi}(x)$ by
 $$
   \sum_{h=0}^{D-1}\chi(h)\frac{te^{t(h+x)}}{e^{Dt}-1}=\sum_{k=0}^\infty
   \frac{B_{k,\chi}(x)}{k!}t^k,
 $$
 that is,
 $$
   B_{k,\chi}(x)=D^{k-1}\sum_{h=0}^{D-1}\chi(h)B_k((h+x)/D)
   =\sum_{j=0}^k\binom kjB_{j,\chi}x^{k-j}.
 $$
 In paricular, $B_{k,\chi}(0)$ is the usual generalized Bernoulli
 number.

 For positive integers $a,c,k,\ell$ satisfying $ka+\ell c=D$ and
 $(a,c)=1$, we choose integers $b$ and $d$ such that $ad-bc=1$ and set
 $$
   \chi(a,c,k,\ell)=\chi(kb+\ell d).
 $$
 (It is easy to see that the definition does not depend on the choice
 of $b$ and $d$.)

 Finally, we let
 $$
   \tau(\chi):=\sum_{h=0}^{D-1}\chi(h)e^{2\pi ih/D}
 $$
 denote the Gaussian sum associated to $\chi$.
 \end{Definition}

   Now we can describe our first main result.

 \begin{Theorem} \label{theorem 1} Let $R_n$, $0<n<w$, be the unique
   cusp form of weight $w+2$ on $\Gamma_0(N)$ characterized by
   $r_n(f)=\gen{f,R_n}$. Let $\chi$ be a primitive Dirichlet character
   modulo $D$ with $D>1$. Then we have
   \begin{equation*}
   \begin{split}
    &r_{\chi}(R_n)(X)+(-1)^{n-1}\chi(-1)r_{\chi}(R_n)(-X) \\
    &\qquad=\frac{(2i)^{w+1}}{\tau(\overline\chi)}\Bigg(
     \epsilon_1(-D)^{-\tn}\frac{B_{\tn+1,\overline\chi}(DX)}{\tn+1}
     -D^{-n}\frac{B_{n+1,\overline\chi}(DX)}{n+1} \\
    &\qquad\qquad+\epsilon_2(-1)^{n-1}\chi(-N)N^\tn
     D^nX^w\frac{B_{\tn+1,\chi}(-1/DNX)}{\tn+1} \\
    &\qquad\qquad+\epsilon_3\chi(-1)D^\tn
     X^w\frac{B_{n+1,\chi}(-1/DX)}{n+1} \\
    &\qquad\qquad
     +G_n(X)+(-1)^{n-1}\chi(-1)G_n(-X)\Bigg),
   \end{split}
   \end{equation*}
   where
   $$
     \epsilon_1=\begin{cases}1, &\text{if }N=1, \\
     0, &\text{if }N>1, \end{cases} \qquad
     \epsilon_2=\begin{cases}1, &\text{if }(N,D)=1, \\
     0, &\text{if }(N,D)>1, \end{cases} \qquad
     \epsilon_3=\begin{cases}1, &\text{if }N|D, \\
     0, &\text{if }N\nmid D, \end{cases}
   $$
   and
   $$
     G_n(X)=\sum_{\substack{a,c,k,\ell>0,(a,c)=1, \\ N|c,ka+\ell c=D}}
     \overline\chi(a,c,k,\ell)(aX+\ell/D)^n(-cX+k/D)^\tn.
   $$
 \end{Theorem}

 In terms of $L$-values, Theorem \ref{theorem 1} can be rephrased as
 follows. Here for a cusp form $f$ in $S_{w+2}(\Gamma_0(N))$, we set
 $$
   \Lambda(f,s)=\left(\frac{2\pi}{\sqrt N}\right)^{-s}\Gamma(s)
   L(f,s), \qquad
   \Lambda(f,\chi,s)=\left(\frac{2\pi}{D\sqrt N}\right)^{-s}
   \Gamma(s)L(f,\chi,s).
 $$
 Note that when $(N,D)=1$, the function $\Lambda(f,\chi,s)$ satisfies a
 functional equation $$\Lambda(f,\chi,s)=\epsilon\Lambda(f,\chi,w+2-s)$$
 for some root of unity $\epsilon$. When $(N,D)>1$, we need to modify
 the definition of $\Lambda$ (replacing $D\sqrt N$ in the denominator
 by other number) in order to get a functional equation of
 the same symmetry. Here we stick to our definition of $\Lambda$ to
 keep the statement of the result simple.

 \begin{Theorem} \label{theorem 2}
   Let $\{f_1,\ldots,f_s\}$ be an orthogonal basis for
   $S_{w+2}(\Gamma_0(N))$, and let $\chi$ be a primitive Dirichlet
   character modulo $D$ with $D>1$. Let $m$ and $n$ be
   integers satisfying $0\le m\le w$, $0<n<w$, and
   $(-1)^{m+n+1}\chi(-1)=1$. Then we have
   \begin{equation*}
   \begin{split}
    &\sum_{i=1}^s\frac1{\gen{f_i,f_i}}\Lambda(f_i,\chi,m+1)
     \overline{\Lambda(f_i,n+1)}
    =(-D)^{m+1}(i\sqrt N)^{m+n+2}r_{m,\chi}(R_n) \\
    &\qquad=\frac{(2i)^{w+1}i^{m+n+2}DN^{(n+1)/2}}{2\binom
      wm\tau(\overline\chi)}
     \times\Bigg(\epsilon_1(-1)^{n+1}\binom\tn\tm D^n
     \frac{B_{\tn-\tm+1,\overline\chi}}{\tn-\tm+1} \\
    &\qquad\qquad
    +\binom n\tm D^\tn\frac{B_{n-\tm+1,\overline\chi}}{n-\tm+1}
    +\epsilon_2(-1)^{n+m}\binom\tn m\chi(-N)N^{\tn-m}D^n
     \frac{B_{\tn-m+1,\chi}}{\tn-m+1} \\
    &\qquad\qquad+\epsilon_3(-1)^{m+1}\binom nm\chi(-1) D^\tn
     \frac{B_{n-m+1,\chi}}{n-m+1} \\
    &\qquad\qquad+2(-1)^{m+1}\sum_{\substack{a,c,k,\ell>0,(a,c)=1, \\
     N|c,D=ka+\ell c}}\overline\chi(a,c,k,\ell) \\
    &\qquad\qquad\qquad\times
     \sum_{r=0}^\tm(-1)^r\binom nr\binom\tn{\tm-r}
     a^rc^{\tm-r}\ell^{n-r}k^{\tn-\tm+r}\Bigg),
   \end{split}
   \end{equation*}
   where $\epsilon_1$, $\epsilon_2$, and $\epsilon_3$ are given as in
   Theorem \ref{theorem 1}.
 \end{Theorem}

 \end{section}

 \begin{section}{Examples}
 \begin{subsection}{Example 1} \label{example: 1}
   Let $N=1$, $D=3$, $w+2=12$, and $\chi=\left(\frac
     {-3}\cdot\right)$. Since $S_{12}(\SL(2,\Z))$ is one-dimensional,
   we expect that for odd $n$ with $1\le n\le 9$,
   the polynomials $g_n(X)=r_\chi(R_n)(X)+(-1)^nr_\chi(R_n)(-X)$ should be
   scalar multiples of each other. Indeed, we have $B_{k,\chi}=0$ for
   even $k$, and
   $$
     B_{1,\chi}=-\frac13, \quad B_{3,\chi}=\frac23, \quad
     B_{5,\chi}=-\frac{10}3, \quad B_{7,\chi}=\frac{98}3, \quad
     B_{9,\chi}=-\frac{1618}3,
   $$
   which give
   \begin{align*}
     B_{2,\chi}(x)&=-\frac23x, \quad B_{4,\chi}(x)=-\frac43x^3+\frac83x,
     \quad B_{6,\chi}(x)=-2x^5+\frac{40}3x^3-20x, \\
     B_{8,\chi}(x)&=-\frac83x^7+\frac{112}3x^5-\frac{560}3x^3
      +\frac{784}3x, \\
     B_{10,\chi}(x)&=-\frac{10}3x^9+80x^7-840x^5+3920x^3-
       \frac{1618}3x.
   \end{align*}
   The tuples $(a,c,k,\ell)$ contributing to $G_n(x)$ are
   $$
    (1,1,1,2),\ (1,1,2,1),\ (1,2,1,1),\ (2,1,1,1)
   $$
   with $\chi(a,c,k,\ell)$ being $-1,1,1,-1$, respectively. We find
   that
   \begin{gather*}
     g_1=-\frac{2048}{\sqrt 3}
     \left(-1536X^9+128X^7-\frac{128}{81}X^3+\frac{512}{2187}X\right), \\
     g_3=g_7=-\frac{25}{48}g_1, \quad g_5=\frac5{12}g_1, \quad g_9=g_1.
   \end{gather*}
   This gives us
   \begin{equation*}
   \begin{split}
     \frac{\Lambda(\Delta,\chi,2)\Lambda(\Delta,2)}{\|\Delta\|^2}
   &=-\frac{2^{18}3^2}5\sqrt3
    =-\frac{\Lambda(\Delta,\chi,10)\Lambda(\Delta,2)}{\|\Delta\|^2} \\
     \frac{\Lambda(\Delta,\chi,4)\Lambda(\Delta,2)}{\|\Delta\|^2}
   &=-\frac{2^{14}3^2}5\sqrt3
    =-\frac{\Lambda(\Delta,\chi,8)\Lambda(\Delta,2)}{\|\Delta\|^2}.
   \end{split}
   \end{equation*}
   and $\Lambda(\Delta,\chi,6)=0$. (Note that the sign for the
   functional equation of $\Lambda(\Delta,\chi,s)$ is $-1$, so that it
   vanishes at $s=6$.)

   The result can be verified numerically as follows. The $L$-values
   can be approximated by the standard method. We have
   \begin{equation*}
   \begin{split}
          \Lambda(\Delta,2)&\approx 0.003707710464948 \\
     \Lambda(\Delta,\chi,2)&\approx -228.22304046813742 \\
     \Lambda(\Delta,\chi,4)&\approx -14.263940029258589.
   \end{split}
   \end{equation*}
   To get an approximation for the inner product $\gen{\Delta,\Delta}$,
   we consider the Poincar\'e series
   $$
     P_k(z)=\frac12\sum_{c,d\in\Z,(c,d)=1}\frac{e^{2\pi
         ik(az+b)/(cz+d)}}{(cz+d)^{12}},
   $$
   where in the summand $a$ and $b$ are any integers satisfying
   $ad-bc=1$. The Poincar\'e series is characterized by the property
   that if $f(z)=\sum_{k=1}^\infty a_k(f)e^{2\pi ikz}$, then
   $$
     \gen{f,P_k}=\frac{\Gamma(11)}{(4\pi k)^{11}}a_f(k).
   $$
   From this we can easily deduce that
   $$
     P_k(z)=2\pi\Gamma(11)\frac{\tau(k)}{(2k)^{11}}\frac{\Delta(z)}
     {\|\Delta\|^2}
   $$
   and
   $$
     \gen{P_k,P_m}=4\pi^2\Gamma(11)^2\frac{\tau(k)\tau(m)}{(4km)^{11}
     \|\Delta\|^2},
   $$
   where $\tau(k)$ is the $k$th Fourier coefficient of $\Delta(z)$.
   Now there is a well-known formula for the inner product $\gen{P_k,P_m}$
   in terms of the Kloosterman sums and the Bessel functions.
   (See \cite[Corollary 3.4]{Iwaniec}.) For instance, evaluating
   $\gen{P_1,P_1}$, we get
   \begin{equation} \label{equation: Delta norm}
     \frac1{\|\Delta\|^2}\approx 965845.709168185.
   \end{equation}
   Then
   $$
     \frac{\Lambda(\Delta,\chi,2)\Lambda(\Delta,2)}{\|\Delta\|^2}
     \approx-817284.10841880
     \approx-\frac{2^{18}3^2}5\sqrt3.
   $$
   Note that the approximation \eqref{equation: Delta norm} can also be
   obtained using the formula
   $$
     \sum_{m=1}^\infty\frac{\tau(m)^2}{m^{20}}=\frac2{245}
     \frac{4^{20}\pi^{29}}{20!}\frac{\zeta(9)}{\zeta(18)}\|\Delta\|^2
   $$
   given in \cite[page 2]{Zagier-VI}.
 \end{subsection}

 \begin{subsection}{Example 2} Let $N=1$ and $w+2=24$. The normalized Hecke
   eigenforms are
 $$
   E_4(z)^3\Delta(z)+(-156\pm 12\sqrt{144169})\Delta(z)^2.
 $$
 Let $f_1$ and $f_2$ denote these two functions. In this example we
 will work out the ratio $\Lambda(f_i,\chi,m)/\Lambda(f_i,12)$ for
 $\chi=\left(\frac5\cdot\right)$ and even $m$. We first express $f_i$
 in terms of $R_n$.

 Let $r_m(f)$ denote $\int_0^{i\infty}f(z)z^m\,dz$. By
 computing the determinant of the matrix
 $$
   \begin{pmatrix}r_1(R_2) & r_3(R_2)\\
   r_1(R_4) & r_3(R_4) \end{pmatrix}
 $$
 using the formula in Theorem 1 of \cite{KZ2}, we easily see that $R_2$
 and $R_4$ form a basis for $S_{24}(\SL(2,\Z))$. Let $T_2$ denote the
 second Hecke operator on $S_{24}(\SL(2,\Z))$. We can determine the
 numbers $a,b,c,d$ such that
 $$
   \begin{pmatrix}T_2R_2\\ T_2R_4\end{pmatrix}
  =\begin{pmatrix}a&b\\c&d\end{pmatrix}\begin{pmatrix}
   R_2\\R_4\end{pmatrix}
 $$
 by considering the relation
 $$
   \begin{pmatrix}a&b\\c&d\end{pmatrix}
  =\begin{pmatrix}r_1(T_2R_2)&r_3(T_2R_2) \\
   r_1(T_2R_4) & r_3(T_2R_4)\end{pmatrix}
   \begin{pmatrix} r_1(R_2) & r_3(R_2) \\
   r_1(R_4) & r_3(R_4) \end{pmatrix}^{-1},
 $$
 which, using the formulas in Theorem 2.8 of \cite{F5}, is shown to be
 $$
   \begin{pmatrix}-716424 &-6894720\\
   1416492/19 & 717504 \end{pmatrix}.
 $$
 From this we deduce that
 $$
   118041R_2+(1135193\pm 19\sqrt{144169})R_4
 $$
 are (unnormalized) Hecke eigenforms, i.e., scalar multiples of $f_1$
 and $f_2$, respectively. Now we have
 $$
   \Lambda(f,\chi,m+1)=(-5i)^{m+1}r_{m,\chi}(f), \qquad
   \Lambda(f,m+1)=(-i)^{m+1}r_m(f)
 $$
 for any cusp form $f$ of weight $24$. Thus, using the formulas from
 Theorem 1 of \cite{KZ2} and our Theorem \ref{theorem 1}, we find that
 \begin{equation*}
 \begin{split}
   \frac{\Lambda(f_i,\chi,2)}{\Lambda(f_i,12)}
 &=\frac{454494815973561283200\mp495053625411273600\sqrt{144169}}
   {11\sqrt5}, \\
   \frac{\Lambda(f_i,\chi,4)}{\Lambda(f_i,12)}
 &=\frac{1710371411434851840\mp1874940923128320\sqrt{144169}}{11\sqrt5},
 \\
   \frac{\Lambda(f_i,\chi,6)}{\Lambda(f_i,12)}
 &=\frac{7923984224047200\mp8900924205600\sqrt{144169}}{11\sqrt 5}, \\
   \frac{\Lambda(f_i,\chi,8)}{\Lambda(f_i,12)}
 &=\frac{46543863219840\mp56895592320\sqrt{144169}}{11\sqrt 5}, \\
   \frac{\Lambda(f_i,\chi,10)}{\Lambda(f_i,12)}
 &=\frac{359949679200\mp545421600\sqrt{144169}}{11\sqrt 5}, \\
   \frac{\Lambda(f_i,\chi,12)}{\Lambda(f_i,12)}
  &=\frac{469261440\mp789120\sqrt{144169}}{\sqrt 5}.
 \end{split}
 \end{equation*}
 Now recall that Theorem 1 of \cite{KZ1} implies that the ratio
 $\Lambda(f_i,\chi,12)/\sqrt 5\Lambda(f_i,12)$ is a square in the ring
 of integers of $\Q(\sqrt{144169})$. Indeed, we find that
 $$
   \frac{\Lambda(f_i,\chi,12)}{\sqrt 5\Lambda(f_i,12)}
  =(3288\mp24\sqrt{144169})^2.
 $$
 \end{subsection}

 \begin{subsection}{Hecke eigenforms on $\SL(2,\Z)$}
 For the convenience of the reader, here we tabulate the (unnormalized)
 Hecke eigenforms in terms of $R_n$ for the case $\dim
 S_k(\SL(2,\Z))=2$. For each weight $k$, we give two bases, one with
 even $n$ and the other with odd $n$.
 $$ \extrarowheight3pt
 \begin{array}{c||l} \hline\hline
 \text{weight} & \text{bases} \\ \hline\hline
 24 & 133705R_1+(1421844\pm12\sqrt{144169})R_3\text{ or } \\
    & 118041R_2+(1135193\pm19\sqrt{144169})R_4 \\ \hline
 28 & 357271915R_1+(5430899304\pm26568\sqrt{18209})R_3\text{ or } \\
    & 166985R_2+(2335719\pm23\sqrt{18209})R_4 \\ \hline
 30 & 339215569R_1+(6031600980\pm6360\sqrt{51349}R_3\text{ or } \\
    & 39282705R_2+(646717136\pm1352\sqrt{51349}R_4 \\ \hline
 32 & 18559684975R_1+(381717886692\pm12876\sqrt{18295489})R_3\text{ or } \\
    & 20837993R_2+(398996469\pm27\sqrt{18295489})R_4 \\ \hline
 34 & 17696951272R_1+(416907865575\pm20925\sqrt{2356201})R_3\text{ or } \\
    & 8056833785R_2+(177566376094\pm17806\sqrt{2356201})R_4 \\ \hline
 38 & 67449635297R_1+(2033146500360\pm\sqrt{63737521})R_3\text{ or } \\
    & 1231612816525R_2+(35003462442636\pm146676\sqrt{63737521})R_4 \\
 \hline\hline
 \end{array}
 $$
 \end{subsection}

 \begin{subsection}{Example 3} Let $N=2$. To obtain exact values of
   ratios between twisted $L$-values of newforms, we can follow the
   following procedure.

   Theorem 1.4 of \cite{FY} asserts that if we let $d_w$ denote the
   dimension of $S_{w+2}(\Gamma_0(2))$, then each of the sets
   $$
     \{R_{2i}:i=1,\ldots,d_w\}, \quad
     \{R_{2i-1}:i=1,\ldots,d_w\}
   $$
   is a basis for $S_{w+2}(\Gamma_0(2))$. Using Theorem 1.1 and Theorem
   1.3 of \cite{FY}, we can find the matrices for the Hecke operators
   with respect to the above bases. Diagonalizing the matrices, we
   obtain expressions of newforms in terms of $R_n$. Then an
   application of Theorem \ref{theorem 1} gives us the values of
   ratios between twisted $L$-values of newforms.

   Let us consider the case $w+2=16$. The space
   $S_{16}(\Gamma_0(2))$ has dimension $3$ and is spanned by
   $$
     f_1=\Delta(z)E_4(z), \quad f_2=\Delta(2z)E_4(2z), \quad
     f_3=\eta(z)^{16}\eta(2z)^{16}.
   $$
   By a direct computation, we find that the unique normalized newform
   is
   $$
     f=f_1+256f_2-600f_3=q-128q^2+6252q^3+\cdots
   $$
   whose eigenvalue for the Atkin-Lehner
   involution $w_2$ is $+1$. Let $D>1$ be a fundamental discriminant. We
   now compute $\Lambda(f,\chi_D,8)/\Lambda(f,8)$ with
   $\chi_D=\left(\frac D\cdot\right)$ for the first few $D$. Note that
   if $(D,2)=1$, then the functional equation for $L(f,\chi_D,s)$ has
   sign $\chi_D(-2)$. Thus, if $D\equiv 5\mod 8$ and $D>0$, we know
   that $L(f,\chi_D,8)=0$.

   Proceeding as in Example 2 and using Theorem 1.1 and Theorem 1.3 of
   \cite{FY}, we find that
   \begin{equation*}
   \begin{split}
     T_3\begin{pmatrix}R_2\\R_4\\R_6\end{pmatrix}
     &=[r_{2j-1}(T_3R_{2i})][r_{2j-1}(R_{2i})]^{-1}
       \begin{pmatrix}R_2\\R_4\\R_6\end{pmatrix} \\
     &=\begin{pmatrix}154348&2478080&3784704 \\
       -11648&-186388&-279552 \\
       1456 & 22880 & 31596 \end{pmatrix}
       \begin{pmatrix}R_2\\R_4\\R_6\end{pmatrix},
 %    &=\begin{pmatrix}r_1(T_3R_2)&r_3(T_3R_2)&r_5(T_3R_2)\\
 %      r_1(T_3R_4)&r_3(T_3R_4)&r_5(T_3R_4) \\
 %      r_1(T_3R_6)&r_3(T_3R_6)&r_5(T_3R_6) \end{pmatrix}
 %      \begin{pmatrix} r_1(R_2)&r_3(R_2)&r_5(R_2) \\
 %      r_1(R_4)&r_3(R_4)&r_5(R_4) \\
 %      r_1(R_6)&r_3(R_6)&r_5(R_6) \end{pmatrix}
   \end{split}
   \end{equation*}
   where $T_3$ denote the third Hecke operator. The characteristic
   polynomial of the above matrix is $(x+3348)^2(x-6252)$. The
   eigenfunction $7R_2+110R_4+168R_6$ associated to the eigenvalue
   $6252$ must be a newform. Applying Theorem 1.1 of
   \cite{FY} and Theorem \ref{theorem 2} we obtain
   $$ \extrarowheight3pt
   \begin{array}{c|c} \hline\hline
   D  & D^{-1/2}\Lambda(f,\chi_D,8)/\Lambda(f,8) \\ \hline
   8  & 2(2^7\cdot3^2)^2 \\
   12 & 2(2^8\cdot3\cdot7)^2 \\
   17 & (2^6\cdot3^2\cdot7^2)^2 \\
   24 & 2(2^8\cdot3\cdot7\cdot29)^2 \\
   28 & 2(2^{10}\cdot3\cdot7\cdot19)^2 \\
   33 & (2^6\cdot3\cdot7\cdot11\cdot23)^2 \\
   40 & 2(2^8\cdot3\cdot5\cdot7\cdot61)^2 \\
   41 & (2^9\cdot3^2\cdot7\cdot23)^2 \\
   44 & 2(2^8\cdot3^2\cdot11\cdot41)^2 \\
   56 & 2(2^9\cdot3^2\cdot5\cdot7^2)^2 \\
   57 & (2^6\cdot3\cdot15671)^2 \\
   60 & 2(2^{10}\cdot3\cdot5\cdot43)^2 \\
   65 & (2^8\cdot3^2\cdot5\cdot13\cdot23)^2 \\
 %  73 & \\
 %  76 & 
   \hline\hline
   \end{array}
   $$
   To check the correctness, we note that a half-integral weight cusp
   form on $\Gamma_0(8)$ corresponding to the normalized newform of
   weight $16$ on $\Gamma_0(2)$ is
   \begin{equation*}
   \begin{split}
     &-\frac{11}{252}[E_6(4\tau),\theta(\tau)]_1+\frac{32}{252}
      [E_6(8\tau),\theta]_1-88\eta(4\tau)^8\eta(8\tau)^8\theta(\tau)
      \\
     &\qquad=q-128q^4-2^7\cdot3^3q^8+4065q^9-2^8\cdot3\cdot7q^{12}
     +2^{14}q^{16}-2^6\cdot3^2\cdot7^2q^{17}+\cdots,
   \end{split}
   \end{equation*}
   where $[g,h]_r$ denotes the Rankin-Cohen bracket, $E_6(\tau)$ is the
   usual Eisenstein series of weight $6$ on $\SL(2,\Z)$, and
   $\theta(\tau)=\sum_{n\in\Z} e^{2\pi in^2\tau}$ is the Jacobi theta
   function.
 \end{subsection}

 \begin{subsection}{Newforms on $\Gamma_0(2)$} For the convenience of
   the reader, here we tabulate newforms on $\Gamma_0(2)$ in terms of
   $R_n$ for the first few $w$. Note that for $w+2=8,10$, we have $\dim
   S_{w+2}(\Gamma_0(2))=1$ and each $R_n$ is a newform. Also, for
   $w+2=12$, the space of newforms has dimension $0$.
   $$ \extrarowheight3pt
   \begin{array}{c||l} \hline\hline
   \text{weight} & \text{bases} \\ \hline\hline
   14 & 21R_1+220R_3,\ R_1+12R_3 \text{ or} \\
      & R_2+8R_4,\ 11R_2+120R_4 \\ \hline
   16 & 49R_1+936R_3+1872R_5 \text{ or}\\
      & 7R_2+110R_4+168R_6 \\ \hline
   18 & 11R_1+300R_3+1056R_5 \text{ or}\\
      & 15R_2+364R_4+1232R_6 \\ \hline
   20 & 11R_1+416R_3+2576R_5+2816R_7, \\
      & \qquad
        3861R_1+123488R_3+321776R_5-622336R_7 \text{ or}\\
      & 51R_2+1722R_4+9464R_6+8448R_8, \ R_2+30R_4+104R_6 \\ \hline
   22 & 143R_1+6612R_3+48640R_5+63232R_7, \\
      & \qquad 1105R_1+52524R_3+425472R_5+708864R_7\text{ or} \\
      & 113R_2+4624R_4+28896R_6+29952R_8, \\
      & \qquad 19R_2+816R_4+6048R_6+9984R_8 \\ \hline
   24 & 10R_3+459R_5+3264R_7+3536R_9 \text{ or} \\
      & 693R_2+34010R_4+228480R_6+16320R_8-311168R_{10} \\ \hline\hline
   \end{array}
   $$

 \begin{Remark} For $N=2$, it is shown in \cite{FY} that
   the first few $R_{2i}$ or the first few $R_{2i-1}$ is a basis for
   $S_{w+2}(\Gamma_0(N))$. Our computation suggests that the same is
   true for $N=3,4,5$. In these cases, we can compute the ratios of
   twisted $L$-values of newforms using the above approach. For $N\ge
   6$, the method no longer works, as the dimension of
   $S_{w+2}(\Gamma_0(N))$ already exceeds the number of $R_n$.
\end{Remark}
\end{subsection}
\end{section}

\begin{section}{Proof of Theorems}
\begin{subsection}{Preliminary} Let $f\in S_{w+2}(\Gamma_0(N))$ and
  $\chi$ be a Dirichlet character modulo $D$ with $D>1$. Recall
  that if $\chi$ is primitive, then we have
  $$
    \chi(n)=\frac1{\tau(\overline\chi)}\sum_{h=0}^{D-1}
    \overline\chi(h)e^{2\pi ihn/D}.
  $$
  It follows that
  $$
    f_\chi(z)=\frac1{\tau(\overline\chi)}\sum_{h=0}^{D-1}
    \overline\chi(h)f(z+h/D)
  $$
  and
  $$
    r_{m,\chi}(f)=\frac1{\tau(\overline\chi)}\sum_{h=0}^{D-1}
    \overline\chi(h)r_{m,h/D}(f),
  $$
  where
  $$
    r_{m,h/D}(f)=\int_0^{i\infty}f(z+h/D)z^m\,dz.
  $$

  Before we proceed to evaluate $r_{m,h/D}(R_n)$, let us recall the
  following properties of the Bernoulli polynomials. (See
  \cite[pages 804--805]{Handbook}.)

\begin{Lemma} \label{lemma: Fourier for Bernoulli} For two real numbers
  $a$ and $x$ and an integer $k$, we have
  $$
    B_k(a+x)=\sum_{j=0}^k\binom kjB_j(a)x^{k-j}.
  $$
  Moreover, the Fourier expansion for the Bernoulli function
  $B_k(\{x\})$, $k\ge 2$, is given by
  $$
    B_k(\{x\})=-\frac{k!}{(2\pi i)^k}\sum_{r\in\Z,r\neq 0}
    \frac{e^{2\pi irx}}{r^k},
  $$
  and
  $$
   -\frac1{2\pi i}\lim_{T\to\infty}\sum_{0<|r|<T}
    \frac{e^{2\pi irx}}{r}
   =\begin{cases}B_1(\{x\})=\{x\}-1/2, &\text{if }x\not\in\Z, \\
    0, &\text{if }x\in\Z. \end{cases}
  $$
\end{Lemma}

  To evaluate $r_{m,h/D}(R_n)$, we shall utilize the following
  expression for $R_n$.

\begin{Lemma} \label{lemma: Rn} Let $R_n(z)$ be the cusp form of
  weight $w+2$ on $\Gamma_0(N)$ characterized by the property
  $r_n(f)=\gen{f,R_n}$ for all $f\in S_{w+2}(\Gamma_0(N))$. We have
$$
  R_n(z)=c_n^{-1}\sum_{\SMat abcd\in\Gamma_0(N)}
  \frac1{(az+b)^{\tn+1}(cz+d)^{n+1}}, \quad
  c_n=(-1)^n4\pi i(2i)^{-w-1}\binom wn.
$$
\end{Lemma}

\begin{proof} See \cite[Lemma 2.1]{FY}.
  (See also \cite[Proposition 1]{AN1}.)
\end{proof}

From the above lemma, we have
$$
  c_nR_n(z+h/D)=\sum_{\SMat abcd\in\Gamma_0(N)}
  \frac1{(az+ah/D+b)^{\tn+1}(cz+ch/D+d)^{n+1}}.
$$
We shall consider the cases
\begin{enumerate}
\item $a=0$, (with $N=1$),
\item $c=0$,
\item $(a,b)=\pm(D,-h)$ (with $(N,D)=1$),
\item $(c,d)=\pm(D,-h)$ (with $N|D$),
\item $ac(ah/D+b)(ch/D+d)<0$,
\item $ac(ah/D+b)(ch/D+d)>0$,
\end{enumerate}
separately. For $j=1,\ldots,6$ and $h\in\Z$ with $(h,D)=1$, we let
$S_{j,h}(z)$ denote the sum of $1/(az+b)^{\tn+1}(cz+d)^{n+1}$ over all
matrices $\SMat abcd$ in $\Gamma_0(N)$ satisfying the $j$th condition.
Note that we have
\begin{equation} \label{equation: Sh=Sh+D}
  S_{j,h+D}(z)=S_{j,h}(z).
\end{equation}
This is because a matrix $\SMat abcd$ contributes to the sum
$S_{j,h+D}(z)$ if and only if the matrix $\SMat a{a+b}c{c+d}$
contributes to the sum $S_{j,h}(z)$. In the following sections, we
will obtain formulas for
$$
  I_{j,h,m}:=\int_0^{i\infty}(S_{j,h}(z)+(-1)^{m+n+1}S_{j,-h}(z))z^m
  \,dz
$$
and
\begin{equation*}
\begin{split}
  F_{j,h}(X)&:=\int_0^{i\infty}S_{j,h}(z)(X-z)^w\,dz
   +(-1)^{n-1}\int_0^{i\infty}S_{j,-h}(z)(X+z)^w\,dz \\
 &=\sum_{m=0}^w(-1)^m\binom wm I_{j,h,m}X^\tm.
\end{split}
\end{equation*}
\end{subsection}

\begin{subsection}{Case $a=0$ (and $N=1$)} In this section, we shall
evaluate the integral
$$
  \int_0^{i\infty}S_{1,h}(z)z^m\,dz.
$$
\begin{subsubsection}{Case $m+1>n$} \label{subsubsection: a=0, m+1>n}
If $a=0$, then
$$
  \Mat abcd=\pm\Mat0{-1}1d, \quad d\in\Z.
$$
We have
\begin{equation*}
\begin{split}
 S_{1,h}(z)
 &=\sum_{\SMat0{-1}1d\in\SL(2,\Z)}\frac{(-1)^{\tn+1}}{(z+h/D+d)^{n+1}}
   +\sum_{\SMat01{-1}d\in\SL(2,\Z)}\frac{(-1)^{n+1}}{(z+h/D-d)^{n+1}} \\
 &=\frac{2(-1)^{n+1}(-2\pi i)^{n+1}}{\Gamma(n+1)}
  \sum_{r=1}^\infty r^ne^{2\pi ir(z+h/D)}.
\end{split}
\end{equation*}
(Here we have used the formula \cite[page 51]{Schoeneberg} for
$\sum_{d\in\Z}(\tau+d)^{-n}$.) It follows that
\begin{equation*}
\begin{split}
  \int_0^{i\infty}S_{1,h}(z)z^m\,dz
% &=(-1)^{n+1}\frac{(-2\pi i)^{n+1}}{n!}i^{m+1}\int_0^\infty
%   S_1(it+h/D)t^m\,dt \\
 &=2(-1)^{n+1}(-2\pi i)^{n-m}\frac{m!}{n!}
   \sum_{r=1}^\infty\frac{e^{2\pi irh/D}}{r^{m-n+1}}.
\end{split}
\end{equation*}
From this and Lemma \ref{lemma: Fourier for Bernoulli}, we obtain
\begin{equation*}
\begin{split}
%  \int^{i\infty}_0(S_{1,h}(z)+(-1)^{m+n+1}S_{1,-h}(z))z^m\,dz
 I_{1,h,m}=\frac{(-1)^m(4\pi i)m!}{n!(m-n+1)!}B_{m-n+1}(\{h/D\})
\end{split}
\end{equation*}
for positive integers $m$ and $n$ with $m+1>n$. (Note that the
integral-sum is no longer absolutely convergent in the case $m=n$.
However, the conclusion remains valid in view of the bounded
convergence theorem.)
\end{subsubsection}

\begin{subsubsection}{Case $m+1<n$} \label{subsubsection: a=0, m+1<n}
It is easy to check that
$$
  S_{1,h}(z)\ll
  \begin{cases} 1, &\text{if }|z|\ll 1, \\
  |z|^{-n},&\text{if }|z|\gg 1. \end{cases}
$$
Therefore, if $m+1<n$, we may integrate term by term. Now we have
\begin{equation} \label{equation: Sh, S-h}
  S_{j,-h}(z)=(-1)^{n+1}S_{j,h}(-z)
\end{equation}
for $j=1,\ldots,6$. Hence,
$$
  \int_0^{i\infty}(S_{1,h}(z)+(-1)^{m+n+1}S_{1,-h}(z))z^m\,dz
 =\int_{-i\infty}^{i\infty}S_{1,h}(z)z^m\,dz.
$$
We then integrate term by term. By shifting the path of integration to
$\operatorname{Re}z=\infty$ or $\operatorname{Re}z=-\infty$
depending on whether $h/D+d$ is positive or negative, we find that the
integral for each term is zero. This shows that the contribution of
$S_{1,h}$ in the case $m+1<n$ is $0$.
\end{subsubsection}

\begin{subsubsection}{Case $m=n-1$} \label{subsubsection: a=0, m=n-1}
From \eqref{equation: Sh, S-h}, we obtain
\begin{equation*}
\begin{split}
 I_{1,h,n-1}=\int_{-i\infty}^{i\infty}S_{1,h}(z)z^{n-1}\,dz
  =2(-1)^{n+1}\lim_{U\to\infty}\sum_{d\in\Z}\int_{-iU}^{iU}\frac{z^{n-1}dz}
  {(z+h/D+d)^{n+1}}.
\end{split}
\end{equation*}
If $n=1$, then we have
\begin{equation*}
\begin{split}
 I_{1,h,n-1}&=-2\lim_{U\to\infty}\sum_{d\in\Z}\left(
    \frac1{iU+h/D+d}-\frac1{-iU+h/D+d}\right) \\
  &=4i\lim_{U\to\infty}\sum_{d\in\Z}\frac U{(h/D+d)^2+U^2}
   =4i\lim_{U\to\infty}\frac1U\sum_{d\in\Z}\frac1
    {((h/D+d)/U)^2+1}.
\end{split}
\end{equation*}
Interpreting the last sum as a Riemann sum, we arrive at
$$
  I_{1,h,n-1}=4i\int_{-\infty}^\infty\frac{dx}{x^2+1}=4\pi i.
$$
If $n\ge 2$, we apply integration by parts once and get
\begin{equation*}
\begin{split}
  I_{1,h,n-1}&=2(-1)^{n+1}\lim_{U\to\infty}\sum_{d\in\Z}\Bigg(
  -\frac{(iU)^{n-1}}{n(iU+h/D+d)^n}
  +\frac{(-iU)^{n-1}}{n(-iU+h/D+d)^n} \\
  &\qquad\qquad\qquad\qquad\qquad\qquad+\frac{n-1}n
   \int_{-iU}^{iU}\frac{z^{n-2}\,dz}{(z+h/D+d)^n}\Bigg).
\end{split}
\end{equation*}
Again, the sum
$$
  \sum_{d\in\Z}\left(
  -\frac{(iU)^{n-1}}{(iU+h/D+d)^n}
  +\frac{(-iU)^{n-1}}{(-iU+h/D+d)^n}\right)
$$
can be interpreted as a Riemann sum for some integral whose
value turns out to be $0$. Thus, we have
$$
  I_{1,h,n-1}=2(-1)^{n+1}\frac{n-1}n\lim_{U\to\infty}\sum_{d\in\Z}
    \int_{-iU}^{iU}\frac{z^{n-2}\,dz}{(z+h/D+d)^n}.
$$
Integrating by parts repeatedly, we eventually obtain
$$
  I_{1,h,n-1}=(-1)^{n+1}\frac{4\pi i}n.
$$
\end{subsubsection}

\begin{subsubsection}{Summary for the case $a=0$} We now combine the
  computations in Sections \ref{subsubsection: a=0, m+1>n},
  \ref{subsubsection: a=0, m+1<n}, and \ref{subsubsection: a=0,
    m=n-1}. We have
\begin{equation*} \label{equation: a=0, temp 1}
\begin{split}
 F_{1,h}(X)&:=\int_0^{i\infty}S_{1,h}(z)(X-z)^w\,dz+(-1)^{n-1}
  \int_0^{i\infty}S_{1,-h}(z)(X+z)^w\,dz \\
&=\sum_{m=0}^w(-1)^m\binom wm I_{1,h,m}X^{w-m}.
\end{split}
\end{equation*}
The result in Section \ref{subsubsection: a=0, m+1>n} shows that
the contribution from the terms with $m\ge n$ is
\begin{equation*}
\begin{split}
 &\frac{4\pi i}{n!}\sum_{m=n}^w\binom wm
  \frac{m!}{(m-n+1)!}X^{w-m}B_{m-n+1}(\{h/D\}) \\
 &\qquad\qquad=\frac{4\pi i}{\tn+1}\binom wn\sum_{m=n}^w
  \binom{w-n+1}{m-n+1}X^{w-m}B_{m-n+1}(\{h/D\}).
\end{split}
\end{equation*}
In view of Lemma \ref{lemma: Fourier for Bernoulli}, this is equal to
$$
  \frac{4\pi i}{\tn+1}\binom wn\left(
  B_{\tn+1}(\{h/D\}+X)-X^{w-n+1}\right).
$$
From Section \ref{subsubsection: a=0, m+1<n}, we know that the
contribution from the terms with $m+1<n$ to \eqref{equation: a=0, temp
  1} is $0$, while Section \ref{subsubsection: a=0, m=n-1} shows
that the term $m=n-1$ yields
$$
  (-1)^{n-1}\binom w{n-1}X^{w-n+1}(-1)^{n+1}\frac{4\pi i}n
 =\frac{4\pi i}{\tn+1}\binom wnX^{w-n+1}.
$$
Combining everything, we get the following formula for $S_{1,h}(z)$.

\begin{Lemma} \label{lemma: a=0} Let $c_n$ be defined as in Lemma
  \ref{lemma: Rn}. We have
$$
  c_n^{-1}I_{1,h,m}
 =(-1)^{n+m}(2i)^{w+1}\binom wm^{-1}\binom\tn\tm
  \frac{B_{\tn-\tm+1}(\{h/D\})}{\tn-\tm+1},
$$
or equivalently,
\begin{equation*}
\begin{split}
  c_n^{-1}F_{1,h}(X)=(-1)^n(2i)^{w+1}\frac{B_{\tn+1}(\{h/D\}+X)}{\tn+1}.
\end{split}
\end{equation*}
\end{Lemma}
\end{subsubsection}
\end{subsection}

\begin{subsection}{Case $c=0$} The evaluation for the case $c=0$ is
  very similar to the case $a=0$, so the proof will be very sketchy.
  We have
  $$
    S_{2,h}(z)=2\sum_{b\in\Z}\frac1{(z+h/D+b)^{\tn+1}}
   =\frac{2(-2\pi i)^{\tn+1}}{\tn!}\sum_{r=1}^\infty r^\tn
    e^{2\pi ir(z+h/D)}.
  $$
  Thus, when $m+1>\tn$,
  $$
    \int_0^{i\infty}S_{2,h}(z)z^m\,dz
   =2(-2\pi i)^{\tn-m}\frac{m!}{\tn!}\sum_{r=1}^\infty
    \frac{e^{2\pi irh/D}}{r^{m-\tn+1}},
  $$
  and
  $$
    I_{2,h,m}
   =\frac{(-1)^{m+n+1}(4\pi i)m!}{\tn!(m-\tn+1)!}B_{m-\tn+1}(\{h/D\}).
  $$
  When $m+1<\tn$, we find that
  $$
    I_{2,h,m}=0
  $$
  The case $m=\tn-1$ yields
  $$
    I_{2,h,\tn-1}=\frac{4\pi i}\tn.
  $$
  In summary, the contribution of the case $c=0$ to the period
  polynomial is the following.

\begin{Lemma} \label{lemma: c=0} We have
$$
  c_n^{-1}I_{2,h,m}=(-1)^{m-1}(2i)^{w+1}\binom wm^{-1}\binom n\tm
  \frac{B_{n-\tm+1}(\{h/D\})}{n-\tm+1},
$$
and
\begin{equation*}
\begin{split}
  c_n^{-1}F_{2,h}(X)=-(2i)^{w+1}\frac{B_{n+1}(\{h/D\}+X)}{n+1}.
\end{split}
\end{equation*}
\end{Lemma}
\end{subsection}

\begin{subsection}{Case $(a,b)=\pm(D,-h)$ (with $(N,D)=1$)} Let $c$
  and $d$ be any two 
  integers satisfying $dD+cNh=1$. Then the matrices in $\Gamma_0(N)$ with
  an upper row $\pm(D,-h)$ are
  $$
    \pm\Mat D{-h}{N(c+kD)}{d-khN}
  $$
  for integers $k\in\Z$. Then
  $$
    S_{3,h}(z)=2\sum_{k\in\Z}\frac1{(Dz)^{\tn+1}((c+kD)Nz+1/D)^{n+1}}.
  $$
  Now we make a change of variable $z\mapsto-1/ND^2z$ in the
  integral
  \begin{equation*}
  \begin{split}
    \int_0^{i\infty}S_{3,h}(z)z^m\,dz&=2
    \int_0^{i\infty}\left(\sum_{k\in\Z}\frac{z^m}
    {(Dz)^{\tn+1}((c+kD)Nz+1/D)^{n+1}}\right)dz \\
   &=2(-1)^{m+n}N^{\tn-m}D^{w-2m}\int_0^{i\infty}
    \left(\sum_{k\in\Z}\frac{z^\tm}{(z-k-c/D)^{n+1}}\right)dz.
  \end{split}
  \end{equation*}
  At this point, we are basically back to the previous cases. For
  $\tm+1>n$, we have
  \begin{equation*}
  \begin{split}
    \int_0^{i\infty}S_{3,h}(z)z^m\,dz
   =2(-1)^{m+n}N^{\tn-m}D^{w-2m}(-2\pi i)^{n-\tm}
    \frac{\tm!}{n!}\sum_{r=1}^\infty\frac{e^{-2\pi irc/D}}{r^{\tm-n+1}}.
  \end{split}
  \end{equation*}
  Then Lemma \ref{lemma: Fourier for Bernoulli} yields
  \begin{equation*}
  \begin{split}
   I_{3,h,m}=-\frac{4\pi iN^{\tn-m}D^{w-2m}\tm!}{n!(\tm-n+1)!}
    B_{\tm-n+1}(\{-c/D\}).
  \end{split}
  \end{equation*}
  For $m$ with $\tm+1<n$, arguing as in Section \ref{subsubsection:
    a=0, m+1<n}, we find that
  $$
    I_{3,h,m}=0.
  $$
  When $\tm+1=n$ (i.e., $m=\tn+1$), a discussion similar to that in Section
  \ref{subsubsection: a=0, m=n-1} leads to
  $$
    I_{3,h,\tn+1}=-\frac{4\pi iD^{w-2\tn-2}}{Nn}.
  $$
  Finally, the condition $dD+cNh=1$ means that $c$ is the
  multiplicative inverse of $Nh$ modulo $D$.

  \begin{Lemma} \label{lemma: (a,b)=(D,-h)} The contribution from
    $S_{3,h}$ is
  $$
    c_n^{-1}I_{3,h,m}=(-1)^{n-1}(2i)^{w+1}\binom wm^{-1}\binom\tn m
    D^{w-2m}N^{\tn-m}\frac{B_{\tn-m+1}(\{-\overline N\overline h/D\})}
    {\tn-m+1}
  $$
  and
  \begin{equation*}
  \begin{split}
   c_n^{-1}F_{3,h}(X)=(-1)^{n-1}(2i)^{w+1}N^\tn D^wX^w\frac{B_{\tn+1}
    (\{-\overline N\overline h/D\}-1/D^2NX)}{\tn+1},
  \end{split}
  \end{equation*}
  where $\overline N$ and $\overline h$ denote the multiplicative
  inverses of $N$ and $h$ modulo $D$, respectively.
  \end{Lemma}
\end{subsection}

\begin{subsection}{Case $(c,d)=\pm(D,-h)$ (with $N|D$)} This case is
  similar to the previous case. Choose integers $a$ and $b$ with
  $ah+bD=-1$. Then the matrices in $\SL(2,\Z)$ with a lower row
  $\pm(D,-h)$ are
  $$
    \pm\Mat{a+kD}{b-kh}D{-h}.
  $$
  Thus,
  $$
    S_{4,h}(\tau)=2\sum_{k\in\Z}\frac1
    {((a+kD)\tau-1/D)^{\tn+1}(D\tau)^{n+1}}.
  $$
  Arguing similarly as in the previous section, we obtain the
  following evaluation.

  \begin{Lemma} \label{lemma: (c,d)=(D,-h)} We have
  $$
    c_n^{-1}I_{4,h,m}=(2i)^{w+1}\binom wm^{-1}\binom nm D^{w-2m}
    \frac{B_{n-m+1}(\{-\overline h/D\})}{n-m+1}
  $$
  and
  \begin{equation*}
  \begin{split}
    c_n^{-1}F_{4,h}=(2i)^{w+1}D^wX^w\frac{B_{n+1}
    (\{-\overline h/D\}-1/D^2X)}{n+1},
  \end{split}
  \end{equation*}
  where $\overline h$ is the multiplicative inverse of $h$ modulo
  $D$.
  \end{Lemma}
\end{subsection}

\begin{subsection}{Case $ac(ah/D+b)(ch/D+d)<0$} From \eqref{equation:
    Sh, S-h} we know that
  \begin{equation} \label{equation: S5h 1}
  \begin{split}
   &\int_0^{i\infty}S_{5,h}(z)(X-z)^w\,dz+(-1)^{n-1}\int_0^{i\infty}
    S_{5,-h}(z)(X+z)^w\,dz \\
   &\qquad\qquad=\int_{-i\infty}^{i\infty}S_{5,h}(z)(X-z)^w\,dz.
  \end{split}
  \end{equation}
  We now assume $h>0$ and evaluate the integral above.

  There are two cases
  $ac>0$ and $ac<0$. In the former case, because $ad-bc=1$, we must
  have $-d/c<h/D<-b/a$. In the latter case we have $-b/a<h/D<-d/c$
  instead. Also, if $\SMat abcd$ contributes to $S_{5,h}$,
  then so does $\SMat{-a}{-b}{-c}{-d}$. It follows that
  \begin{equation*}
  \begin{split}
    S_{5,h}(z)&=2\sum_{\substack{\SMat abcd\in\Gamma_0(N),a,c>0,\\
    -d/c<h/D<-b/a}}\frac1{(az+ah/D+b)^{\tn+1}(cz+ch/D+d)^{n+1}} \\
   &\qquad\qquad+2\sum_{\substack{\SMat abcd\in\Gamma_0(N),a>0,c<0, \\
    -b/a<h/D<-d/c}}\frac1{(az+ah/D+b)^{\tn+1}(cz+ch/D+d)^{n+1}}.
  \end{split}
  \end{equation*}
  Now the condition $ad-bc=1$ implies that $-b/a$ and $-d/c$
  are Farey neighbors. Then a general property of the Farey
  fractions says that in order for a fraction $h/D>0$ to be sandwiched
  between $-b/a$ and $-d/c$, $h$ and $D$ must be of the form
  $h=k|b|+\ell|d|$ and $D=k|a|+\ell|c|$ for some positive integers $k$
  and $\ell$. This in particular shows that the number of terms in the
  sum $S_{5,h}$ is finite.

  In the case $a,c>0$, the integers $b$ and $d$ are non-positive.
  Thus, $D=ka+\ell c$ and $h=-kb-\ell d$. We have
  $$
    \frac{ah}D+b=\frac{a(-kb-\ell d)}{ka+\ell c}+b=-\frac\ell{ka+\ell c}
    =-\frac\ell D.
  $$
  Likewise, we have
  $$
    \frac{ch}D+d=\frac kD.
  $$
  In the case $a>0$, $c<0$, we have $b\le 0$ and $d\ge 0$. Thus,
  $D=ka-\ell c$, $h=-kb+\ell d$, and
  $$
    \frac{ah}D+b=\frac\ell D, \qquad
    \frac{ch}D+d=\frac kD.
  $$
  Then $S_{5,h}$ becomes
  \begin{equation} \label{equation: S5h}
  \begin{split}
    S_{5,h}(z)&=2\sum_{\substack{\SMat abcd\in\Gamma_0(N),a,c,k,\ell>0,\\
    D=ka+\ell c,h=-kb-\ell d}}\frac1{(az-\ell/D)^{\tn+1}(cz+k/D)^{n+1}} \\
   &\qquad\qquad+2\sum_{\substack{\SMat abcd\in\Gamma_0(N),a,k,\ell>0,c<0, \\
    D=ka-\ell c,h=-kb+\ell d}}\frac1{(az+\ell D)^{\tn+1}(cz+k/D)^{n+1}}.
  \end{split}
  \end{equation}
  Now let us recall a formula from \cite{FY}.

  \begin{Lemma} Let $a,b,c,d$ be real numbers such that $abcd<0$. Then
  $$
    \int_{-i\infty}^{i\infty}\frac{(X-\tau)^w}{(a\tau+b)^{\tn+1}
    (c\tau+d)^{n+1}}\,d\tau
   =\frac{(-1)^n2\pi i}{(ad-bc)^{w+1}}\binom wn
    \mathrm{sgn}(ab)(aX+b)^n(cX+d)^\tn.
  $$
  \end{Lemma}

  \begin{proof} See \cite[page 341]{FY}.
  \end{proof}

  Applying this lemma to \eqref{equation: S5h}, we obtain the
  following formula for \eqref{equation: S5h 1}.

  \begin{Lemma} \label{lemma: <0} We have
  \begin{equation*}
  \begin{split}
    c_n^{-1}F_{5,h}(X)=&-(2i)^{w+1}
     \sum_{\substack{\SMat abcd\in\Gamma_0(N),a,c,k,\ell>0,\\
       D=ka+\ell c,h=-kb-\ell d}}(aX-\ell/D)^n(cX+k/D)^\tn \\
    &+(2i)^{w+1}
     \sum_{\substack{\SMat abcd\in\Gamma_0(N),a,k,\ell>0,c<0, \\
      D=ka-\ell c,h=-kb+\ell d}}(aX+\ell/D)^n(cX+k/D)^\tn.
  \end{split}
  \end{equation*}
  \end{Lemma}
\end{subsection}

\begin{subsection}{Case $ac(ah/D+b)(ch/D+d)>0$} The evaluation of the
  terms with $ac(ah/D+b)(ch/D+d)>0$ follows the argument in
  \cite[page 13--16]{FY}. Here we only provide a sketch.

  Firstly, we have $S_{6,-h}(z)=(-1)^{n-1}S_{6,h}(-z)$, so that
  $$
    I_{6,h,m}=\int_{-i\infty}^{i\infty}S_{6,h}(z)z^m\,dz
  $$
  We can show that
  \begin{equation*}
  \begin{split}
   &\sum_{ac(ah/D+b)(ch/D+d)>0}
    \frac1{|(az+ah/D+b)^{\tn+1}(cz+ch/D+d)^{n+1}|} \\
   &\qquad\qquad\qquad
    \ll\begin{cases} 1/|z|, &\text{if }|z|\ll 1, \\
    1/|z|^{w+1}, &\text{if }|z|\gg 1.\end{cases}
  \end{split}
  \end{equation*}

  \begin{subsubsection}{Case $0<m<w$} If $m$ is not $0$ or $w$, we may
    change the order of integration of summation. In this case, since
    the two poles of $1/(az+ah/D+b)^{\tn+1}(cz+ch/D+d)^{n+1}$ lie on
    the same side of the imaginary axis, we have
  $$
    I_{6,h,m}=0.
  $$
  \end{subsubsection}

  \begin{subsubsection}{Case $m=w$} We have
  \begin{equation*}
  \begin{split}
  I_{6,h,w}&=\lim_{\epsilon\to 0} \sum_{ac(ah/D+b)(ch/D+d)>0}
    \int_{-i/\epsilon}^{i\epsilon}\frac{z^w\,dz}
    {(az+ah/D+b)^{\tn+1}(cz+ch/D+d)^{n+1}} \\
   &=-\lim_{\epsilon\to 0}\sum\left(\int_{i/\epsilon}^{i\infty}+
    \int_{-i\infty}^{-i/\epsilon}\right)
   (\text{same function})\,dz.
  \end{split}
  \end{equation*}
  Making a change of variable $z\mapsto i/(\epsilon t)$, we obtain
  $$
    I_{6,h,w}=i\lim_{\epsilon\to 0}\epsilon
      \sum\int_{-1}^1\frac{dt}{(a-i\epsilon t(ah/D+b))^{\tn+1}
      (c-i\epsilon t(ch/D+d))^{n+1}}.
  $$
  For a given pair of integers $a$ and $c$, we fix integers $b_0$ and
  $d_0$ such that $ad_0-b_0c=1$. The other integers $b$ and $d$
  satisfying $ad-bc=1$ are $b=b_0+ak$ and $d=d_0+ck$ for $k\in\Z$.
  Then
  \begin{equation*}
  \begin{split}
  I_{6,h,w}&=i\lim_{\epsilon\to 0}\epsilon\sum_{a,c}
    \frac1{a^{\tn+1}c^{n+1}}\sum_{u\in\Z+b_0/a+h/D}\int_{-1}^1
    \frac{dt}{(1-i\epsilon tu)^{\tn+1}(1-i\epsilon t(u+1/ac))^{n+1}}
    \\
  &=i\sum_{a,c}\frac1{a^{\tn+1}c^{n+1}}
    \int_{-\infty}^{\infty}\int_{-1}^1\frac{dt}{(1-itu)^{w+2}}\,du
  \end{split}
  \end{equation*}
  (Note that there might be some integers $k$ such that
  $ac(ah/D+b_0+ak)(ch/D+d_0+ck)$ is not positive. However, it should
  be clear that the contribution of these terms to the above sum is
  not significant.) Following the computation in \cite[page 341]{FY},
  we find that
  $$
    \int_{-\infty}^\infty\int_{-1}^1\frac{dt}{(1-ixt)^{w+2}}\,dx
   =\frac{2\pi}{w+1},
  $$
  and
  $$
    I_{6,h,w}=\frac{2\pi i}{w+1}\sum_{(a,c)=1,N|c,a,c\neq 0}
       \frac1{a^{\tn+1}c^{n+1}}.
  $$
  The sum is equal to $0$ if $n$ is even, and is equal to
  $$
%  &=\frac{8\pi i}{w+1}
%    \sum_{a=1,(a,N)=1}^\infty\frac1{a^{\tn+1}}\sum_{c=1,(a,c)=1}^\infty
%     \frac1{(cN)^{n+1}} \\
  I_{6,h,w}=\frac{8\pi i}{w+1}\frac{\zeta(n+1)\zeta(\tn+1)}
    {\zeta(w+2)N^{n+1}}\prod_{p|N}\frac{1-p^{-(\tn+1)}}{1-p^{-(w+2)}}
  $$
  if $n$ is odd. In either case, we have
  \begin{equation} \label{equation: S6 1}
    I_{6,h,w}=-\frac{4\pi i}{N^{n+1}}\binom
    wn\frac{w+2}{B_{w+2}}\frac{B_{n+1}}{n+1}\frac{B_{\tn+1}}{\tn+1}
    \prod_{p|N}\frac{1-p^{-\tn-1}}{1-p^{-w-2}},
  \end{equation}
  where $p$ runs through all prime divisors of $N$.
  This settles the case $m=w$.
  \end{subsubsection}

  \begin{subsubsection}{Case $m=0$} We first make a change of variable
    $z\to-1/ND^2z$ and obtain
  \begin{equation*}
  \begin{split}
   I_{6,h,0}=\int_{-i\infty}^{i\infty}\sum\frac{(-1)^nN^\tn D^wz^w\,dz}
    {(a/D-N(ah+bD)z)^{\tn+1}(-c/ND+(ch+dD)z)^{n+1}}.
  \end{split}
  \end{equation*}
  Now choose integers $u$ and $v$ such that $Du-Nhv=1$.
  For each $\left(\begin{smallmatrix}a&b\\c&d\end{smallmatrix}\right)$
  in $\Gamma_0(N)$, we set
  $$
    \begin{pmatrix}\alpha&\beta\\\gamma&\delta\end{pmatrix}
   =\begin{pmatrix}a&b\\c&d\end{pmatrix}
    \begin{pmatrix}u&h\\Nv&D\end{pmatrix}.
  $$
  We check that the condition $ac(ah/D+b)(ch/D+d)>0$ holds if and only
  if the matrix
  $$
    \begin{pmatrix}a'&b'\\c'&d'\end{pmatrix}
   :=\begin{pmatrix}\delta&-\gamma/N\\-\beta N&\alpha\end{pmatrix}
  $$
  satisfies
  $$
    a'c'(a'v/D+b')(c'v/D+d')>0.
  $$
  In fact, we have
  $$
    a'=ch+dD, \quad c'=-N(ah+bD), \quad
    a'v/D+b'=-c/DN, \quad c'v/D+d'=a/D,
  $$
  and as $\left(\begin{smallmatrix}
  a&b\\c&d\end{smallmatrix}\right)$ goes through every element in
  $\Gamma_0(N)$ satisfying $ac(ah/D+b)(ch/D+d)>0$, the corresponding
  $\left(\begin{smallmatrix}a'&b'\\c'&d'\end{smallmatrix}\right)$
  goes through elements in $\Gamma_0(N)$ satisfying
  $a'c'(a'v/D+b')(c'v/D+d')>0$. It follows that
  \begin{equation*}
  \begin{split}
   I_{6,h,0}
   =\int_{-i\infty}^{i\infty}\sum_{a',b',c',d'}\frac{(-1)^nN^\tn
     D^wz^w\,dz}{(a'z+a'v/D+b')^{n+1}(c'z+c'v/D+d')^{\tn+1}}.
  \end{split}
  \end{equation*}
  By \eqref{equation: S6 1}, this is equal to
  $$
    \frac{4\pi iD^w}N\binom wn
    \frac{w+2}{B_{w+2}}\frac{B_{n+1}}{n+1}\frac{B_{\tn+1}}{\tn+1}
    \prod_{p|N}\frac{1-p^{-n-1}}{1-p^{-w-2}}.
  $$

  \begin{Remark} Note that the argument above can be extended to show
  that 
  $$
    r_{m,h/D}(R_n)=(-1)^{n+m}N^{\tn-m}D^{\tm-m}r_{\tm,v/D}
    (R_\tn)
  $$
  for $0\le m\le w$ and $1\le n\le w-1$. Here the integer $v$ is the
  multiplicative inverse of $-Nh$ modulo $D$ since $u$ and $v$ satisfy
  $Du-Nhv=1$.
  \end{Remark}
  \end{subsubsection}

  \begin{subsubsection}{Summary for the case $ac(ah/D+b)(ch/D+d)>0$}
  Combining the computations above, we arrive at the following
  conclusion.

  \begin{Lemma} \label{lemma: S6} We have
  \begin{equation*}
  \begin{split}
   c_n^{-1}F_{6,h}(X)&=
    (-1)^n(2i)^{w+1}\frac{w+2}{B_{w+2}}\frac{B_{n+1}}{n+1}
    \frac{B_{\tn+1}}{\tn+1} \\
   &\qquad\qquad\times
    \left(X^w\frac{D^w}N\prod_{p|N}\frac{1-p^{-n-1}}{1-p^{-w-2}}
   -\frac1{N^{n+1}}\prod_{p|N}\frac{1-p^{-\tn-1}}{1-p^{-w-2}}\right),
  \end{split}
  \end{equation*}
  where the products run over all prime divisors $p$ of $N$.
  \end{Lemma}
  \end{subsubsection}
\end{subsection}
\end{section}

\begin{subsection}{Proof of Theorem \ref{theorem 1}} This is just a
  summarization of Lemmas \ref{lemma: a=0}, \ref{lemma: c=0},
  \ref{lemma: (a,b)=(D,-h)}, \ref{lemma: (c,d)=(D,-h)},
  \ref{lemma: <0}, and \ref{lemma: S6}.
\end{subsection}

\begin{subsection}{Proof of Theorem \ref{theorem 2}}
  Since $\{f_1,\ldots,f_s\}$ is an orthogonal basis, we have
  \begin{equation} \label{equation: Rn in basis}
    R_n=\sum_{i=1}^s\frac{\gen{R_n,f_i}}{\gen{f_i,f_i}}f_i.
  \end{equation}
  Now applying $r_{m,\chi}$ to both sides of \eqref{equation: Rn in
    basis}, we obtain
  \begin{equation*}
  \begin{split}
    r_{m,\chi}(R_n)&=\sum_{i=1}^s\frac{\gen{R_n,f_i}}{\gen{f_i,f_i}}
      r_{m,\chi}(f_i)
     =\sum_{i=1}^s\frac1{\gen{f_i,f_i}}\overline{r_n(f_i)}
      r_{m,\chi}(f_i) \\
    &=\frac{n!}{(2\pi i)^{n+1}}\frac{m!}{(-2\pi i)^{m+1}}
      \sum_{i=1}^s\frac1{\gen{f_i,f_i}}
      \overline{L(f_i,n+1)}L(f_i,\chi,m+1).
  \end{split}
  \end{equation*}
  Then Theorem \ref{theorem 2} follows from Lemmas \ref{lemma: a=0},
  \ref{lemma: c=0}, \ref{lemma: (a,b)=(D,-h)},
  \ref{lemma: (c,d)=(D,-h)}, \ref{lemma: <0}, and \ref{lemma: S6}.
\end{subsection}

\end{document}